\documentclass[a4paper,11pt]{article}
\usepackage{amsmath}
\usepackage{amsfonts}
\usepackage{amsthm}
\usepackage{amssymb}
\usepackage{latexsym}
\usepackage[latin1]{inputenc}
\usepackage{makeidx}

\newcommand{\bbt}{\mathbb{T}}
\newcommand{\bbd}{\mathbb{D}}

\newcommand{\bbz}{\mathbb{Z}}

\newcommand{\calc}{\mathcal{C}}

\newcommand{\call}{\mathcal{L}}
\newcommand{\dsp}{\displaystyle}
\newcommand{\e}{\varepsilon}

\newtheorem{theo}{Theorem}[section]

\newtheorem{cor}[theo]{Corollary}
\newtheorem{lemme}[theo]{Lemma}

\begin{document}

\title{On the growth of powers of operators with spectrum contained in Cantor sets}
\author{AGRAFEUIL Cyril}
\date{}
\maketitle

\begin{abstract}
For $\xi \in \big( 0, \frac{1}{2} \big)$, we denote by $E_{\xi}$ the perfect symmetric set associated to $\xi$, that is
$$
E_{\xi} = \Big\{ \exp \big( 2i \pi (1-\xi) \dsp \sum_{n = 1}^{+\infty} \epsilon_{n} \xi^{n-1} \big) : \, \epsilon_{n} = 0 \textrm{ or } 1 \quad (n \geq 1) \Big\}.
$$
Let $s$ be a nonnegative real number, and $T$ be an invertible bounded operator on a Banach space with spectrum included in $E_{\xi}$. We show that if
\begin{eqnarray*} 
& & \big\| T^{n} \big\| = O \big( n^{s} \big), \,n \rightarrow +\infty \\
& \textrm{and} & \big\| T^{-n} \big\| = O \big( e^{n^{\beta}} \big), \, n \rightarrow +\infty \textrm{ for some } \beta < \frac{\log{\frac{1}{\xi}} - \log{2}}{2\log{\frac{1}{\xi}} - \log{2}},
\end{eqnarray*}
then for every $\e > 0$, $T$ satisfies the stronger property
$$
\big\| T^{-n} \big\| = O \big( n^{s+\frac{1}{2}+\e} \big), \, n \rightarrow +\infty.
$$  
This result is a particular case of a more general result concerning operators with spectrum satisfying some geometrical conditions. 
\end{abstract}

\section{Introduction}

We denote by $\bbt$ the unit circle and by $\bbd$ the open unit disk.  We shall say that a closed subset $E$ of $\bbt$ is a $K$-set if there exists a positive constant $c$ such that for any arc $L$ of $\bbt$,
$$
\sup_{z \in L} d(z,E) \geq c |L|, \eqno (K)
$$
where $|L|$ denotes the length of the arc $L$ and $d(z,E)$ the distance between $z$ and $E$. Let $E$ be a $K$-set. We set
$$
\delta(E) = \sup \Big\{ \delta \geq 0 : \, \int_{0}^{2\pi} \frac{1}{d(e^{it},E)^{\delta}} \mathrm{d}t < +\infty \Big\}.
$$
We have $\dsp \delta(E) \geq \frac{\log{\frac{1}{1-c}}}{\log{\frac{2}{1-c}}}$ (see \cite{Dyn} section 5, proof of lemma 2 and corollary). E. M. Dyn'kin showed in \cite{Dyn} that condition $(K)$ characterizes the interpolating sets for $\Lambda_ {s}^{+}(\bbt)$, $s>0$ (see section 2 for the definition of $\Lambda_ {s}^{+}(\bbt)$). Let $s$ be a nonnegative real number, and let $T$ be an invertible operator on a Banach space. We show (theorem \ref{operateurs}) that if the spectrum of $T$ is included in $E$ and if $T$ satisfies 
\begin{eqnarray*} 
& & \big\| T^{n} \big\| = O \big( n^{s} \big), \,n \rightarrow +\infty \\
& \textrm{and} & \big\| T^{-n} \big\| = O \big( e^{n^{\beta}} \big), \, n \rightarrow +\infty \textrm{ for some } \beta < \frac{\delta(E)}{1+\delta(E)},
\end{eqnarray*}
then for every $\e>0$, $T$ also satisfies the stronger property
\begin{eqnarray} \label{intro1}
\big\| T^{-n} \big\| = O \big( n^{s+\frac{1}{2}+\e} \big), \, n \rightarrow +\infty.
\end{eqnarray}
For $\dsp \xi \in \Big( 0, \frac{1}{2} \Big)$, we denote by $E_{\xi}$ the perfect symmetric set associated to $\xi$, that is
$$
E_{\xi} = \Big\{ \exp \big( 2i \pi (1 - \xi) \dsp \sum_{n = 1}^{+\infty} \epsilon_{n} \xi^{n-1} \big) : \, \epsilon_{n} = 0 \textrm{ or } 1 \quad (n \geq 1) \Big\}.
$$
We set $\dsp b(\xi) = \frac{\log{\frac{1}{\xi}} - \log{2}}{2\log{\frac{1}{\xi}} - \log{2}}$. We obtain (as a consequence of theorem \ref{operateurs}) that if the spectrum of $T$ is included in $E_{\xi}$, $\big\| T^{n} \big\| = O \big( n^{s} \big), \, n \rightarrow +\infty$ and $\big\| T^{-n} \big\| = O \big( e^{n^{\beta}} \big), \, n \rightarrow +\infty$ for some $\beta < b(\xi)$, then $T$ satisfies (\ref{intro1}). Notice that J. Esterle showed in \cite{Est2} that if $T$ is a contraction on a Banach space (respectively on a Hilbert space) with spectrum included in $\dsp E_{\frac{1}{q}}$ (respectively included in $E_{\xi}$) such that $\big\| T^{-n} \big\| = O \big( e^{n^{\beta}} \big), \, n \rightarrow +\infty$ for some $\dsp \beta < b \big( \frac{1}{q} \big)$ (respectively $\beta < b(\xi)$), then $\dsp \sup_{n \geq 0} \big\| T^{-n} \big\| < +\infty$ (respectively $T$ is an isometry). Here $q$ is an integer greater than or equal to $3$. \\

\section{Growth of powers of operators}

Let $p$ be a non-negative integer. We denote by $\calc^{p}(\bbt)$ the space of $p$ times continuously differentiable functions on $\bbt$. We set 
$$
\textrm{{\LARGE $a$}}^{p}(\bbd) = \Big\{ f \in \calc^{p}(\bbt) : \, \widehat{f}(n) = 0 \quad (n<0) \Big\},
$$
$\calc^{\infty}(\bbt)=\bigcap_{p\geq 0} \calc^{p}(\bbt)$ and $\textrm{{\LARGE $a$}}^{\infty}(\bbd)= \bigcap_{p\geq 0} \textrm{{\LARGE $a$}}^{p}(\bbd)$.

\noindent Let $s$ be a nonnegative real number, we denote by $[s]$ the nonnegative integer such that $[s] \leq s < [s]+1$. We define the Banach algebra
$$
\Lambda_{s}(\bbt) = \Big\{ f \in \calc^{[s]}(\bbt) : \, \sup_{z,z' \in \bbt} \frac{\big| f^{([s])}(z) - f^{([s])}(z') \big|}{|z-z'|^{s-[s]}}< +\infty \Big\},
$$
equiped with the norm $\dsp \big\| f \big\|_{\Lambda_{s}} = \big\| f \big\|_{\calc^{[s]}(\bbt)} + \sup_{z,z' \in \bbt} \frac{\big| f^{([s])}(z) - f^{([s])}(z') \big|}{|z-z'|^{s-[s]}}$. We also define the subalgebra
$$
\lambda_{s}(\bbt) = \Big\{ f \in \calc^{[s]}(\bbt) : \, \big| f^{([s])}(z) - f^{([s])}(z') \big| = o \big( |z-z'|^{s-[s]} \big), \, |z-z'| \rightarrow 0 \Big\},
$$
which we equip with the same norm. We also set
\begin{eqnarray*}
\Lambda_{s}^{+}(\bbt) & = & \Big\{ f \in \Lambda_{s}(\bbt) : \, \widehat{f}(n) = 0 \quad (n<0) \Big\} \\
\textrm{and } \lambda_{s}^{+}(\bbt) & = & \Big\{ f \in \lambda_{s}(\bbt) : \, \widehat{f}(n) = 0 \quad (n<0) \Big\}.
\end{eqnarray*}
We remark that if $s$ is an integer, $\Lambda_{s}(\bbt) = \lambda_{s}(\bbt) = \calc^{s}(\bbt)$ and so $\Lambda_{s}^{+}(\bbt) = \lambda_{s}^{+}(\bbt) = \textrm{{\LARGE $a$}}^{s}(\bbd)$. We define 
$$
N_{s}(E) = \big\{ f \in \Lambda_{s}(\bbt) : \, f_{|_{E}} = \ldots = f_{|_{E}}^{([s])} = 0 \big\},
$$
and set $N_{s}^{+}(E) = N_{s}(E) \cap \Lambda_{s}^{+}(\bbt)$. 

\begin{lemme} \label{bernstein}
Let $s$ be a nonnegative real number. Then for all $\e > 0$, we have the following continuous embedding
$$
\Lambda_{s+\frac{1}{2}+\e}(\bbt) \hookrightarrow A_{s}(\bbt).
$$
\end{lemme} 

\begin{proof}
For $s=0$, this is a result of Bernstein (see \cite{Kaha1}, p.13). The general case is obtained by the same arguments. Let $\e > 0$, and set $\tilde{s}=s+\frac{1}{2}+\e$. Let $f \in \Lambda_{\tilde{s}}(\bbt)$. For $h>0$, define
$$
P(h) = \int_{0}^{2\pi} \big| f^{([\tilde{s}])}(e^{i(t-h)}) -  f^{([\tilde{s}])}(e^{i(t+h)}) \big|^{2} \mathrm{d} t.
$$
It follows from Parseval equality that
\begin{eqnarray} \label{parseval}
P(h) = 8 \pi \sum_{n=-\infty}^{+\infty} \big| \widehat{f^{([\tilde{s}])}}(n) \big|^{2} \sin^{2}{(nh)}.
\end{eqnarray}
Let $j_{0}$ be the smallest integer such that $[\tilde{s}] < 2^{j_{0}}$ and let $j \geq j_{0}$. It follows from the relation $\widehat{f^{([\tilde{s}])}}(n) = \Big( \prod\limits_{k=1}^{[\tilde{s}]} (n+k) \Big) \widehat{f}(n+[\tilde{s}]) \, (n \in \bbz)$ and from (\ref{parseval}) that there exists a constant $C_{1} > 0$ independent of $f$ such that 
\begin{eqnarray} \label{bernstein1}
P(h) & \geq & \frac{4}{C_{1}^{2}} \sum_{|n|=2^{j}}^{2^{j+1}-1} \big| \widehat{f}(n+[\tilde{s}]) \big|^{2} (1+|n|)^{2[\tilde{s}]} \sin^{2}{(nh)}.
\end{eqnarray} 
Using the Cauchy-Schwartz inequality, we have
\begin{eqnarray} \label{bernstein2}
\sum_{|n|=2^{j}}^{2^{j+1}-1} \big| \widehat{f}(n+[\tilde{s}]) \big| (1+|n|)^{s} \leq \Big( \sum_{|n|=2^{j}}^{2^{j+1}-1} \big| \widehat{f}(n+[\tilde{s}]) \big|^{2} (1+|n|)^{2[\tilde{s}]} \Big)^{\frac{1}{2}} \Big( \sum_{n=2^{j}}^{2^{j+1}-1} (1+|n|)^{2s-2[\tilde{s}]} \Big)^{\frac{1}{2}}
\end{eqnarray}
Set $\dsp h=\frac{\pi}{3.2^{j}}$. For all integers $n$ such that $2^{j} \leq |n| \leq 2^{j+1}-1$, we have $\dsp \frac{\pi}{3} \leq |nh| \leq \frac{2\pi}{3}$, and so $\dsp \sin^{2}{(nh)} \geq \frac{1}{4}$. So, we deduce from (\ref{bernstein1}) that
$$
\Big( \sum_{|n|=2^{j}}^{2^{j+1}-1} \big| \widehat{f}(n+[\tilde{s}]) \big|^{2} (1+|n|)^{2[\tilde{s}]} \Big)^{\frac{1}{2}} \leq C_{1} P \big( \frac{\pi}{3.2^{j}} \big)^{\frac{1}{2}}.
$$
Then, as $f \in \Lambda_{\tilde{s}}(\bbt)$, we have
$$
P \big( \frac{\pi}{3.2^{j}} \big)^{\frac{1}{2}} \leq (2 \pi)^{\frac{1}{2}} \big\| f \big\|_{\Lambda_{\tilde{s}}} \big( \frac{2\pi}{3.2^{j}} \big)^{\tilde{s}-[\tilde{s}]},
$$
so that
\begin{eqnarray} \label{bernstein3}
\Big( \sum_{|n|=2^{j}}^{2^{j+1}-1} \big| \widehat{f}(n+[\tilde{s}]) \big|^{2} (1+|n|)^{2[\tilde{s}]} \Big)^{\frac{1}{2}} \leq  C_{1} (2 \pi)^{\frac{1}{2}} \big\| f \big\|_{\Lambda_{\tilde{s}}} \big( \frac{2\pi}{3.2^{j}} \big)^{\tilde{s}-[\tilde{s}]}.
\end{eqnarray}
Furthermore, there exists a constant $C_2>0$ such that
\begin{eqnarray} \label{bernstein4} 
\Big( \sum_{|n|=2^{j}}^{2^{j+1}-1} (1+|n|)^{2s-2[\tilde{s}]} \Big)^{\frac{1}{2}} \leq C_2 2^{j \big( s-[\tilde{s}]+\frac{1}{2} \big)}.
\end{eqnarray} 
Finally we deduce from (\ref{bernstein2}) and the inequalities (\ref{bernstein3}) and (\ref{bernstein4}) that there exists a constant $C_{3} > 0$ independent of $f$ such that for all $j \geq j_{0}$,
$$
\sum_{|n|=2^{j}}^{2^{j+1}-1} \big| \widehat{f}(n+[\tilde{s}]) \big| (1+|n|)^{s} \leq 2^{j ( s-\tilde{s}+\frac{1}{2} )} C_{3} \big\| f \big\|_{\Lambda_{\tilde{s}}} = 2^{-\e j} C_{3} \big\| f \big\|_{\Lambda_{\tilde{s}}}.
$$
Summing over $j \geq j_{0}$ these inequalities, we get
$$
\sum_{|n| \geq 2^{j_{0}}} \big| \widehat{f}(n+[\tilde{s}]) \big| (1+|n|)^{s} \leq \frac{C_{3}}{1-2^{-\e}} \big\| f \big\|_{\Lambda_{\tilde{s}}},
$$
On the other hand, we have $\big| \widehat{f}(n) \big| \leq \big\| f \big\|_{\Lambda_{\tilde{s}}}$ for every $n \in \bbz$. So, since $j_{0}$ is independent of $f$, there exists a constant $K > 0$ (independent of $f$) such that
$$
\big\| f \big\|_{s} \leq K \big\| f \big\|_{\Lambda_{\tilde{s}}}
$$ 
\end{proof}

Before giving the main theorem of the paper, we need the following lemma.

\begin{lemme} \label{sansfacteur}
Let $E$  be a closed subset of $\bbt$. We assume that there exists $\delta > 0$ for which $\dsp \int_{0}^{2\pi} \frac{1}{d(e^{it},E)^{\delta}} \mathrm{d}t < +\infty$. Let $\dsp \beta < \frac{\delta}{1+\delta}$ and let $T$ be an invertible operator on a Banach space with spectrum included in $E$ that satisfies 
\begin{eqnarray*} 
& & \big\| T^{n} \big\| = O \big( n^{s} \big), \,n \rightarrow +\infty \quad(\textrm{ for some nonnegative real s}) \\
& \textrm{and} & \big\| T^{-n} \big\| = O \big( e^{n^{\beta}} \big), \, n \rightarrow +\infty ,
\end{eqnarray*}
Then there exists an outer function $f \in \textrm{{\LARGE $a$}}^{\infty}(\bbd)$ which vanishes exactly on $E$ and such that $f(T) := \sum\limits_{n=0}^{+\infty} \widehat{f}(n) T^{n}=0$.
\end{lemme}

\begin{proof}
Let $\omega$ be the weight defined by $\omega(n) = \big\| T^{n} \big\| \, (n \in \bbz)$. Let $\Phi$ be the continuous morphism from $A_{\omega}(\bbt)$ to $\call(X)$ defined by 
$$
\Phi(f) = f(T) = \sum_{n = - \infty}^{+\infty} \widehat{f}(n) T^{n} \qquad \big( f \in A_{\omega}(\bbt) \big).
$$
Since the algebra $A_{\omega}(\bbt)$ is regular, we have $\big\{ z \in \bbt : \, f(z)=0 \quad (f \in \textrm{Ker } \Phi) \big\} \subset E$ (see \cite{EStZo1}, theorem 2.5), and so $J_{\omega}(E) \subset \textrm{Ker } \Phi$. Then the result follows from lemmas 7.1 and 7.2 of \cite{Est2}.
\end{proof}

\begin{theo} \label{operateurs} 
Let $E$ be a $K$-set, and let $s$ be a nonnegative real number. Then, any invertible operator $T$ on a Banach space with spectrum included in $E$ that satisfies 
\begin{eqnarray*} 
& & \big\| T^{n} \big\| = O \big( n^{s} \big), \,n \rightarrow +\infty \\
& \textrm{and} & \big\| T^{-n} \big\| = O \big( e^{n^{\beta}} \big), \, n \rightarrow +\infty \textrm{ for some } \beta < \frac{\delta(E)}{1+\delta(E)},
\end{eqnarray*}
also satisfies the stronger property
$$
\big\| T^{-n} \big\| = O \big( n^{s+\frac{1}{2}+\e} \big), \, n \rightarrow +\infty,
$$  
for all $\e>0$.
\end{theo}

\begin{proof}
Let $\e>0$ and set $\tilde{s}=s+\frac{1}{2}+\e$. Without loss of generality, we may assume that $\tilde{s}$ is not an integer. Let $t$ a real number, which is not an integer, and satisfies $s + \frac{1}{2} < t < \tilde{s}$ and $[t] = [\tilde{s}]$. According to lemma \ref{bernstein}, we can define a continuous morphism $\Phi$ from $\lambda_{t}^{+}(\bbt)$ to $\call(X)$ by 
$$
\Phi(f) = f(T) = \sum_{n=0}^{+\infty} \widehat{f}(n) T^{n} \qquad \big( f \in \lambda_{t}^{+}(\bbt) \big).
$$
Let $I = \textrm{Ker } \Phi$, $I$ is a closed ideal of $\lambda_{t}^{+}(\bbt)$. We denote by $S_{I}$ its inner factor, that is the greatest common divisor of all inner factors of the non-zero functions in $I$ (see \cite{Hoff} p.85), and we set, for $0 \leq k \leq [t]$, $\dsp h^{k}(I) = \big\{ z \in \bbt : \, f(z) = \ldots = f^{(k)}(z) = 0 \quad (f \in I) \big\}$. \\F. A. Shamoyan showed in \cite{Sha} that
$$
I = \Big\{ f \in \lambda_{t}^{+}(\bbt) : \, S(I) | \, S(f) \textrm{ and } f^{(k)}=0 \textrm{ on } h^{k}(I) \textrm{ for all } 0 \leq k \leq [t] \Big\},
$$ 
where $S(f)$ denotes the inner factor of $f$ and $S(I) | \, S(f)$ means that $S(f)/ S(I)$ is a bounded holomorphic function in $\bbd$. Since $\dsp \beta < \frac{\delta(E)}{1+\delta(E)}$, there exists $0 < \delta < \delta(E)$ such that $\dsp \beta < \frac{\delta}{1+\delta}$. We have, by definition of $\delta(E)$, $\dsp \int_{0}^{2\pi} \frac{1}{d(e^{it},E)^{\delta}} \mathrm{d}t < +\infty$. So we deduce from lemma \ref{sansfacteur} that there exists an outer function $f \in \textrm{{\LARGE $a$}}^{\infty}(\bbd)$ which vanishes exactly on $E$ and such that $f \in I$. Therefore, we have $S(I)=1$ and $h^{0}(I) \subset E$, so that $N_{t}^{+}(E) \cap \lambda_{t}(\bbt) \subset I$. Now, as $\Lambda_{\tilde{s}}^{+}(\bbt) \subset \lambda_{t}^{+}(\bbt)$, we can define a continuous morphism $\Psi$ from $\Lambda_{\tilde{s}}^{+}(\bbt)$ to $\call(X)$ by $\Psi = \Phi_{|_{\Lambda_{\tilde{s}}^{+}(\bbt)}}$. Using what precedes, we have 
$$
N_{\tilde{s}}^{+}(E) \subset \textrm{Ker } \Psi.
$$
So there exists a continuous morphism $\tilde{\Psi}$ from $\Lambda_{\tilde{s}}^{+}(\bbt) / N_{\tilde{s}}^{+}(E)$ into $\call(X)$ such that $\Psi = \tilde{\Psi} \circ \pi_{\tilde{s}}^{+}$, where $\pi_{\tilde{s}}^{+}$ is the canonical surjection from $\Lambda_{\tilde{s}}^{+}(\bbt)$ to $\Lambda_{\tilde{s}}^{+}(\bbt) / N_{\tilde{s}}^{+}(E)$. Since $E$ is a $K$-set, by a theorem of E. M. Dyn'kin \cite{Dyn}, it is an interpolating set for $\Lambda_{\tilde{s}}^{+}(\bbt)$, so that the canonical imbedding $i$ from $\Lambda_{\tilde{s}}^{+}(\bbt) / N_{\tilde{s}}^{+}(E)$ into $\Lambda_{\tilde{s}}(\bbt) / N_{\tilde{s}}(E)$ is onto. We have, for $n \geq 0$,
$$
T^{-n} = \tilde{\Psi} \circ i^{-1} \circ \pi_{\tilde{s}} (\alpha^{-n}),
$$
where $\pi_{\tilde{s}}$ denote the canonical surjection from $\Lambda_{\tilde{s}}(\bbt)$ to $\Lambda_{\tilde{s}}(\bbt) / N_{\tilde{s}}(E)$ and where $\alpha \ :\ z \to z$ is the identity map. So we have, for $n \geq 0$,
\begin{eqnarray*}
\big\| T^{-n} \big\| & \leq & \big\| \tilde{\Psi} \circ i^{-1} \big\| \big\|  \pi_{\tilde{s}} (\alpha^{-n}) \big\|_{\Lambda_{\tilde{s}}} \\
& \leq & \big\| \tilde{\Psi} \circ i^{-1} \big\| (1+n)^{\tilde{s}},
\end{eqnarray*}
which completes the proof.
\end{proof}

We give two immediate corollaries of this theorem.

\begin{cor} \label{operateurs1}
Let $\xi \in \big( 0, \frac{1}{2} \big)$ and let $s$ be a nonnegative real number. Then, any invertible operator $T$ on a Banach space with spectrum included in $E_{\xi}$ that satisfies 
\begin{eqnarray*} 
& & \big\| T^{n} \big\| = O \big( n^{s} \big), \,n \rightarrow +\infty \\
& \textrm{and} & \big\| T^{-n} \big\| = O \big( e^{n^{\beta}} \big), \, n \rightarrow +\infty \textrm{ for some } \beta < b(\xi),
\end{eqnarray*}
also satisfies the stronger property
$$
\big\| T^{-n} \big\| = O \big( n^{s+\frac{1}{2}+\e} \big), \, n \rightarrow +\infty,
$$  
for all $\e>0$.
\end{cor}

\begin{proof}
It is well known that $E_{\xi}$ is a $K$-set (see proposition 2.5 of \cite{Est2}). Moreover, $E_{\xi}$ satisfies $\dsp
\int_{0}^{2\pi} \frac{1}{d(e^{it},E)^{\delta}} \mathrm{d}t < +\infty$ if and only if $\dsp \delta < 1+\frac{\log{2}}{\log{\xi}}$. Indeed, the condition $\dsp
\int_{0}^{2\pi} \frac{1}{d(e^{it},E)^{\delta}} \mathrm{d}t < +\infty$ is equivalent to $\sum\limits_{n=1}^{+\infty} \sum\limits_{i=1}^{2^{n-1}} \big| L_{n,i} \big|^{1-\delta}$, where $L_{n,i}$ are the arcs contiguous to $E_{\xi}$, and $\big| L_{n,i} \big|$ are their length, which is equal to $2 \pi \xi^{n-1} (1-2\xi)$ (see \cite{KaSa} for further details). Then it is easily seen that the last series converges if and only if $\dsp \delta < 1+\frac{\log{2}}{\log{\xi}}$, so $\dsp \delta(E_{\xi}) = 1+\frac{\log{2}}{\log{\xi}}$. Now, the result follows immediately from theorem \ref{operateurs}.
\end{proof}

Then we obtain an other immediate result, which generalizes theorem 4.1 of \cite{ElKe}. Indeed, the condition "$\big\| T^{-n} \big\| = O \big( e^{n^{\beta}} \big), \, n \rightarrow +\infty$" which appears in the following corollary is weaker than the condition used by the authors of \cite{ElKe}.

\begin{cor} \label{operateurs2}
Let $E$ be a $K$-set, and let $s$ be a nonnegative real number. Then, there exists a constant $\beta > 0$ independent of $s$ such that any invertible operator $T$ on a Banach space with spectrum included in $E$ that satisfies 
\begin{eqnarray*} 
& & \big\| T^{n} \big\| = O \big( n^{s} \big), \,n \rightarrow +\infty \\
& \textrm{and} & \big\| T^{-n} \big\| = O \big( e^{n^{\beta}} \big), \, n \rightarrow +\infty,
\end{eqnarray*}
also satisfies the stronger property
$$
\big\| T^{-n} \big\| = O \big( n^{s+\frac{1}{2}+\e} \big), \, n \rightarrow +\infty,
$$  
for all $\e>0$.
\end{cor}

\begin{proof}
As $E$ is a $K$-set, we deduce from \cite{Dyn} (section 5, corollary) that $\delta(E)>0$. Then the result follows immediately from theorem \ref{operateurs},  with any $\beta < \frac{\delta(E)}{1+\delta(E)}$.
\end{proof}

\hspace*{-6mm}\textbf{Remark 2.6:} \\
1) Some results concerning operators with countable spectrum are obtained in \cite{Zar2} and in \cite{bibi1}. Let $E$ be a closed subset of $\bbt$ and let $s$, $t$ be two nonnegative reals. We denote by $P(s,t,E)$ the following property: every invertible operator $T$ on a Banach space such that $\textrm{Sp} \, T \subset E$ and satisfies the conditions:
\begin{eqnarray*} 
& & \big\| T^{n} \big\| = O(n^{s}) \,\, (n \rightarrow +\infty) \label{introCp} \\
& & \big\| T^{-n} \big\| = O(e^{\e \sqrt{n}}) \,\, (n \rightarrow +\infty), \, \textrm{ for all } \e > 0, 
\end{eqnarray*}
also satisfies the stronger property
\begin{eqnarray*}
\big\| T^{-n} \big\| = O(n^{t}) \, (n \rightarrow +\infty). 
\end{eqnarray*}
M. Zarrabi showed in \cite{Zar2} (th\'eor\`eme 3.1 and remarque 2.a) that a closed subset $E$ of $\bbt$ satisfies $P(0,0,E)$ if and only if $E$ is countable. Notice that $E$ is called a Carleson set if $\dsp \int_{0}^{2\pi} \log^{+} \frac{1}{d(e^{it}, E)} \mathrm{d}t < +\infty$.
If $E$ is a countable closed subset of $\bbt$, we show in \cite{bibi1} that the following conditions are equivalent: \\
\hspace*{3mm} (i) there exist two positive constants  $C_1, C_2$ such that for every arc  $I\subset \bbt$,
$$
\frac{1}{|I|} \int_{I} \log^{+} \frac{1}{d(e^{it}, E)} \mathrm{d}t \leq C_{1} \log {\frac{1}{|I|}} + C_{2}. \\
$$
\hspace*{2mm} (ii) $E$ is a Carleson set and for all $s \geq 0$, there exists $t$ such that $P(s,t,E)$ is satisfied. \\
For contractions with spectrum satisfying the Carleson condition, we can see \cite{Kell}. \\ 
2) When $\dsp \xi=\frac{1}{q}$, the constant $\dsp b \big( \frac{1}{q} \big)$ in corollary \ref{operateurs1} is the best possible in view of \cite{ERaZa}, where the authors built a contraction $T$ such that $\dsp \lim_{n \rightarrow +\infty} \log \big\| T^{-n} \big\| = +\infty$, $Sp \, T \subset E_{\frac{1}{q}}$ and $\log \big\| T^{-n} \big\| = O \big( n^{b(\frac{1}{q})} \big)$. According to theorem 6.4 of \cite{Est1}, $T$ doesn't satisfy $\big\| T^{-n} \big\| = O \big( n^{s} \big)$ for any real $s \geq 0$. \\

\ \\

\hspace*{-7mm} 2000 Mathematics Subject Classification: 46J15, 46J20, 47A30. \\
Key-words: operators, Beurling algebra, spectral synthesis, perfect symmetric set.

\ \\

\hspace*{-7mm} AGRAFEUIL Cyril \\
Cyril.Agrafeuil@math.u-bordeaux.fr \\
Laboratoire Bordelais d'Analyse et G\'eom\'etrie (LaBAG), CNRS-UMR 5467
Universit\'e Bordeaux I \\
351, cours de la lib\'eration \\
33405 Talence cedex, FRANCE. \\


\begin{thebibliography}{1}


\bibitem{bibi1} C. Agrafeuil, \emph{Id\'eaux ferm\'es de certaines alg\`ebres de Beurling et application aux op\'erateurs \`a spectre d\'enombrable}, preprint.
 
\bibitem{Dyn} E. M. Dynkin, \emph{Free interpolation set for H\"older classes}, Mat. Sbornik, \textbf{109} (1979), 107-128.

\bibitem{ElKe} O. El-Fallah et K. Kellay, \emph{Sous-espaces biinvariants pour certains shifts pond\'er\'es}, Ann. Inst. Fourier \textbf{48} (1998), no. 5, 1543-1558.

\bibitem{Est1} J. Esterle, \emph{Uniqueness, strong form of uniqueness and negative powers of contractions}, Banach Center Publ. \textbf{30} (1994), 1-19.

\bibitem{Est2} J. Esterle, \emph{Distributions on Kronecker sets, strong forms of uniqueness, and closed ideals of $A^{+}$}, J. reine angew. Math. \textbf{450} (1994), 43-82.

\bibitem{ERaZa} J. Esterle, M. Rajoelina and M. Zarrabi, \emph{On contractions with spectrum contained in Cantor set}, Math. Proc. Camb. Phil. Soc. \textbf{117} (1995), 339-343.

\bibitem{EStZo1} J. Esterle, E. Strouse and F. Zouakia, \emph{Theorems of Katznelson-Tzafriri type for contractions}, J. Func. Anal. (2) \textbf{94} (1990), 273-287.

\bibitem{Hoff} K. Hoffman \emph{''Banach spaces of analytic functions''}, Prentic-Hall, Englewood Cliffs, 1962. 

\bibitem{Kaha1} J. P. Kahane, \emph{''S\'eries de Fourier absolument convergentes''}, Erg. Math. \textbf{336}, Springer Verlag, Berlin-Heidelberg-New York, 1973.

\bibitem{KaSa} J. P. Kahane, R. Salem, \emph{''Ensembles parfaits et s\'eries trigonom\'etriques''}, Paris, Hermann, 1963.

\bibitem{Kell} K. Kellay, \emph{Contractions et hyperdistributions \`a spectre de Carleson}, J. London Math. Soc. (2) \textbf{58} (1998), 185-196.

\bibitem{Sha} F. A. Shamoyan \emph{Closed ideals in algebras of functions analytic in the disc and smooth up to its boundary}, Mat. Sbornik \textbf{79} (1994), 425-445.

\bibitem{Zar2} M. Zarrabi, \emph{Contractions \`a spectre d\'enombrable et propri\'et\'e d'unicit\'e des ferm\'es d\'enombrables du cercle unit\'e}, Ann. Inst. Fourier \textbf{43} (1993), 251-263.


\end{thebibliography}
\end{document}